\documentclass[12pt,a4paper,twoside]{article}
\usepackage{amssymb,amsthm}
\usepackage[namelimits,sumlimits,fleqn]{amsmath}
\usepackage{setspace}
\usepackage[a4paper,top=35mm, bottom=30mm, left=20mm, right=20mm]{geometry} %   M A R G I N S
\usepackage[small, margin=20pt]{caption}
\usepackage{float}\restylefloat{figure}
\usepackage{url,hyperref}

\usepackage{fancyhdr}
\allowdisplaybreaks[4]

\title{The ${L^1}$-norm of exponential sums in ${Z^d}$}
\author{Giorgis Petridis}
\date{}

\theoremstyle{plain}
\newtheorem{theorem}{Theorem}[section]
\newtheorem{lemma}[theorem]{Lemma}

\theoremstyle{definition}

\newtheorem{question}[theorem]{Question}

\newtheorem*{definition}{Definition}
\newtheorem*{acknowledgement}{Acknowledgement}

\theoremstyle{definition}
\newtheorem*{remark}{Remark}

\begin{document}

\pagenumbering{arabic}

\setcounter{section}{0}

\bibliographystyle{plain}

\maketitle

\thispagestyle{plain}

\begin{abstract}
Let $A$ be a finite set of integers and $F_A(x) = \sum_{a\in A} \exp(2\pi i a x)$ be its exponential sum. McGehee, Pigno $\&$ Smith and Konyagin have independently proved that $\|F_A\|_1\geq c \log|A|$ for some absolute constant $c$. The lower bound has the correct order of magnitude and was first conjectured by Littlewood. In this paper we present lower bounds on the $L^1$-norm of exponential sums of sets in the $d$-dimensional grid $\mathbb{Z}^d$. We show that $\|F_A\|_1$ is considerably larger than $\log|A|$ when $A\subset \mathbb{Z}^d$ has multidimensional structure. %More precisely we prove results of the following kind. Let $A$ be a random subset of $\{1,\dots, N\}\times\{1,\dots,N\}$, where every element is chosen independently with probability 1/2. Then for all $\varepsilon>0$, $\|F_A\|_1\geq \log^{3/2-\varepsilon} N$. 
We furthermore prove similar lower bounds for sets in $\mathbb{Z}$, which in a technical sense are multidimensional and discuss their connection to an inverse result on the theorem of McGehee, Pigno $\&$ Smith and Konyagin.
\end{abstract}

\section[Introduction]{Introduction}
\label{Introduction}

\vspace{1ex}

We begin with a notational remark. Throughout the paper expressions of the form $Q \leq C$ are taken to mean that the quantity $Q$ is less than an appropriately chosen absolute constant $C>1$. We will therefore write counter-intuitive statements like $2C\leq C$. When the constant is less than 1 a lower case $c$ is used.  

For finite $A\subset \mathbb{Z}^d$ the \textit{exponential sum of} $A$ is $$ F_A(x) = \sum_{a\in A} e(a \cdot x)\;,$$ where $\cdot$ is the usual dot product in $\mathbb{R}^d$, $e(t)=\exp(2\pi i t)$ and $x$ lies in the $d$-dimensional torus $\mathbb{T}^d$. The \textit{$L^1$-norm} of $F_A$ is given by $$ \|F_A\|_1 = \int_{x\in \mathbb{T}^d} |F_A(x)| \, dx\;.$$ We will also write $$ \langle f , g \rangle =   \int_{x\in \mathbb{T}^d} f(x) \overline{g(x)}\, dx$$ for the inner product of two functions $f,g : \mathbb{T}^d \mapsto \mathbb{C}$.

J.E. Littlewood conjectured in 1948 \cite{Littlewood1948} that for all finite sets $A\subset \mathbb{Z}$: $$ \| F_A \|_1 \geq c \log|A|\;.$$ The conjecture was proved in 1980 independently by O.C. McGehee, L. Pigno $\&$ B. Smith \cite{MPS1981} and S.V. Konyagin \cite{Konyagin1981}. 
\begin{theorem}[McGehee--Pigno--Smith, Konyagin] \label{MPSG}
Let $A$ be a finite sets of integers. Then $$\|F_A\|_1\geq c \log|A|\;.$$
\end{theorem} 
Taking $A$ to be a symmetric arithmetic progression about zero, and hence $F_A$ the Dirichlet kernel, shows that the lower bound is of the correct order of magnitude \cite{Lorch1954}. 

The first proof works equally well when $A\subset\mathbb{Z}^d$. The order of magnitude of the lower bound is attained when $A$ is an arithmetic progression in $\mathbb{Z}^d$. On the other hand, if $A$ is the $d$-dimensional cube $\{(x_1,\dots,x_d) : 1\leq x_i \leq N~\mbox{for all $i$}\} \subset \mathbb{Z}^{d}$, then $\|F_A\|_1 = \|F_{\{1,\dots,N\}}\|_1^d \geq (C \log N)^d$. It is therefore natural to ask whether a similar lower bound on $\|F_A\|_1$ holds when $A$ has a genuinely multidimensional structure. 

We answer this question to the affirmative, not only for sets in $\mathbb{Z}^d$, but also for sets in $\mathbb{Z}$. Our results present partial progress towards answering a question of W.T. Gowers on the $L^1$-norm of exponential sums in $\mathbb{Z}^2$, which will be stated below. They also help characterise sets of integers $A$ for which $\|F_A\|_1$ is nearly minimal.

The first step is to quantify what we mean by `genuinely multidimensional structure'. The most typical example that comes to mind is that of the $d$-dimensional cube, where as we have seen $\|F_A\|_1$ is roughly speaking $\log^d|A|$. %It is probable that a similar lower bound holds when $A$ is a random dense subset of the $d$-dimensional cube. 
The identity $\|F_A\|_1 = \|F_{\{1,\dots,N\}}\|_1^d$ no longer holds when $A$ is tweaked and taken to be $\{(a_1+x_1,\dots,a_d+x_d) : 1\leq x_i \leq N~\mbox{for all $i$}\}$ for fixed integers $a_1,\dots,a_d$. We study $\|F_A\|_1$ for sets that have a similar structure and show that in this case $\|F_A\|_1\geq \log^{c d}|A|$. To keep the notation simple, here and most importantly in the proofs that follow, we will from now on set $d=2$ or 3. Our methods can be generalised in a straightforward manner for $d>3$. Considering the general case would make what already is a notation-heavy argument even more technical without adding anything to the method. 

Let us now introduce some terminology, which will be helpful in pinning down an exact meaning for `multidimensional structure'.
\begin{definition}
Let $j\in\{1,2,3\}$,  $a_i\in\mathbb{Z}$ for $i\in \{1,2,3\}\setminus\{j\}$ and $A\subseteq \mathbb{Z}^3$. The intersection of $A$ with the line $\{(x_1,x_2,x_3) : x_i=a_i~\mbox{for $i\in  \{1,2,3\}\setminus\{j\}$}\}$ is a \textit{row of $A$}. 
\end{definition}
\begin{definition}
Let $i\in\{1,2,3\}$, $a_i\in\mathbb{Z}$ and $A \subseteq \mathbb{Z}^3$. The intersection of $A$ with the plane $\{(x_1,x_2,x_3) : x_i=a_i\}$ is a \textit{planar slice of $A$}.
\end{definition}
We call $A\subseteq \mathbb{Z}^2$ a \textit{genuinely 2-dimensional} set, if its rows are either empty or large. We call $A\subseteq \mathbb{Z}^3$ a a \textit{genuinely 3-dimensional} set, if its planar slices are either empty or a genuinely 2-dimensional set.

The first of our results asserts that, if $A$ is genuinely 2-dimensional then $\|F_A\|_1$ is considerably larger than $\log|A|$.  
\begin{theorem} \label{2Dlittlewood}
Let $A\subset \mathbb{Z}^2$ be finite. Suppose that $A$ consists of at least $r$ rows of size at least $s$. Then $$\|F_A\|_1 \geq c \log s \left(\frac{\log r}{\log\log r}\right)^{1/2}.$$
\end{theorem}
The stated lower bound is probably not best possible. Gowers has asked whether $\|F_A\|_1\geq c \log r \log s$ holds. Theorem~\ref{2Dlittlewood} only gives $\|F_A\|_1\geq \log s \log^{1/2-\varepsilon}r$ for all $\varepsilon>0$ and sufficiently large $A$.

%As a corollary to the theorem we get a lower bound on $\|F_A\|_1$ for the case when $A$ is a dense random subset of the 2-dimensional cube.
%\begin{corollary}\label{Random}
%Let $N$ a fixed positive integer and $A$ be a random subset of $\{(x,y) : 1\leq x,y\leq N\}$ where elements are chosen independently with probability $1/2$. Then with high probability $$\|F_A\|_1 \geq c\, \frac{\log^{3/2}N}{\log\log^{1/2}N}\;. $$   
%\end{corollary}

The method of proof of Theorem~\ref{2Dlittlewood} can also be applied to subsets of $\mathbb{Z}$. To define `multidimensional structure' in the integers we turn to a notion often used in additive problems. 
\begin{definition}
Let $A$ and $B$ be sets in two additive groups. A map $$\theta : A\mapsto B$$ is a \textit{Freiman isomorphism of degree $k$} if it is a bijection and $a_1+ \dots + a_k = a_{k+1}+ \dots +a_{2k}$ holds if and only if $\theta(a_1)+ \dots +\theta(a_k) = \theta(a_{k+1})+\dots+ \theta(a_{2k})$ holds for any choice of $a_1,\dots,a_{2k}\in A$. We say $A$ is \textit{Freiman isomorphic of degree $k$} to $B$. 
\end{definition}
Our second main result asserts that if $A\subset \mathbb{Z}$ is Freiman isomorphic to a 3-dimensional set in $\mathbb{Z}^3$, then $\|F_A\|_1$ is considerably larger than $\log|A|$. 
\begin{theorem}\label{Freiman Isomorphism Invariance}
Let $A\subset \mathbb{Z}^3$ be finite. Suppose that  $A$ consists of at least $p$ planar slices each in turn consisting of at least $r$ rows of size at least $s$. If $B\subset \mathbb{Z}$ is Freiman isomorphic of degree $k$ to $A$, then  $$ \| F_B \|_1 \geq c \left(\frac{\log s\,\log r\,\log p}{\log\log s\, \log\log r \,\log\log p} \right)^{1/2}\;,$$ provided that $k= 62 \log r \log s \log p$.
\end{theorem}
A helpful, if imprecise, way to rephrase the above is that $\|F_B\|_1\geq \log^{3/2-\varepsilon}|B|$ for all $\varepsilon>0$ whenever $B\subset \mathbb{Z}$ is isomorphic to a genuinely 3-dimensional set in $\mathbb{Z}^3$ and is sufficiently large. As a consequence we see that any sufficiently large set $A$ where $\|F_A\|_1 \leq C \log|A|$ cannot have this particular 3-dimensional structure. 

The lower bound in Theorem~\ref{Freiman Isomorphism Invariance} is probably not best possible. Moreover, one suspects that the conclusion holds for smaller values of $k$. It is furthermore likely that if $A$ is Freiman isomorphic to a 2-dimensional set in $\mathbb{Z}^2$, then $\|F_A\|_1 \geq \log^{1+\eta}|A|$ for some absolute $0<\eta\leq 1$. The method we present is not strong enough to prove this.

The remaining sections are organised as follows. In Sect.~\ref{Cohen Davenport Pichorides} we prove a lemma that is central to the proof of both theorems. The lemma is a generalisation of a method developed by P.J. Cohen \cite{Cohen1960} to tackle Littlewood's conjecture and was later refined by H. Davenport \cite{Davenport1960} and S.K. Pichorides \cite{Pichorides1974}. In Sect.~\ref{2 dim Littlewood} we prove Theorem~\ref{2Dlittlewood} %and Corollary~\ref{Random}
. In Sect.~\ref{Freiman Isomorphisms} we prove Theorem~\ref{Freiman Isomorphism Invariance}. Finally, in Sect.~\ref{Additive Structure when L1 Norm Small} we discuss how an inverse result for Theorem~\ref{MPSG} may look like and compare the suggested structure with that which comes out of Theorem~\ref{Freiman Isomorphism Invariance}.

\begin{acknowledgement}
The author would like to thank Tim Gowers for proposing the question and for generously sharing his insight. Many of the ideas used in the paper originated from conversations with him. He would also like to thank Ben Green for many helpful discussions and suggestions. 
\end{acknowledgement}

\section[A method of Cohen, Davenport and Pichorides]{A method of Cohen, Davenport and Pichorides}

\label{Cohen Davenport Pichorides}

To prove Theorem~\ref{2Dlittlewood} and Theorem~\ref{Freiman Isomorphism Invariance} we will rely on a combination of techniques developed to tackle Littlewood's conjecture by Cohen \cite{Cohen1960}, Davenport \cite{Davenport1960}, Pichorides \cite{Pichorides1974} and McGehee, Pigno $\&$ Smith \cite{MPS1981}. The four aforementioned papers on the Littlewood conjecture concentrate on constructing a test function $g$ that satisfies two properties: $\|g\|_\infty \leq 1$ and $\langle g, F_A \rangle \geq \log^\alpha |A|$ for some absolute constant $\alpha$. This immediately gives $\log^\alpha |A|\leq\langle g,F_A\rangle \leq \|g\|_\infty \|F_A\|_1\leq \|F_A\|_1$.

Our strategy to prove Theorem~\ref{2Dlittlewood} is as follows. For simplicity let us assume that $A$ consists of $r$ rows $A_1,\dots,A_r$ of size at least $s$, where $A_i\subset \{(x,n_i) : x\in\mathbb{Z}\}$ for some integers $n_1,\dots,n_r$. Let $\Phi_{n_i}$ be the McGehee--Pigno--Smith test function for the exponential sum $F_{A_i}$. That is the function constructed by McGehee, Pigno and Smith that satisfies the two properties listed above for $\alpha=1$. We will combine these to produce a better test function for $A$. This will be done by mirroring the method of Cohen, Davenport and Pichorides.   

Cohen combined the exponentials $\{e(nx): n\in A\}$ and obtained a test function which yields the value $\alpha=1/8-\varepsilon$ for all $\varepsilon>0$. Davenport improved this to $\alpha=1/4-\varepsilon$ and Pichorides to $\alpha=1/2-\varepsilon$. The three arguments are rather similar. A closer look at the underlying method reveals that one can get the same result even when relaxing the most commonly used properties of exponentials to:
\begin{itemize}
\item$ |e(n x)|\leq 1$ for all $n$ and $x$.
\item$\langle e(nx)\,e(mx) , e(kx) \rangle = 0$ unless $k= n+m$. 
\item$\langle e(nx),e(nx)\rangle \geq c$ for all $n$.
\end{itemize}
Our strategy is to replace the exponentials in the existing proofs by the $\Phi_{n_i}$, which satisfy the first condition. The support of the Fourier transform $\widehat{\Phi_{n_i}}$ lies in the line that contains $A_i$ and therefore the $\Phi_{n_i}$ also satisfy the following new versions of the later two conditions.
\begin{itemize}
\item Let $k$ and $l$ be positive integers. $\langle F_A , \Phi_{n_{i_1}} \Phi_{n_{i_2}} \cdots \Phi_{n_{i_k}} \overline{ \Phi_{n_{i_{{k+1}}}} \Phi_{n_{i_{{k+2}}}} \cdots \Phi_{n_{i_{{k+l}}}}} \rangle = 0$
unless $n_{i_1} + \dots + n_{i_k} - n_{i_{k+1}}- \dots -n_{i_{k+l}} = n_\nu$ for some $1\leq \nu\leq R$.
\item$\langle\Phi_{n_i}, F_A \rangle = \langle\Phi_{n_i}, F_{A_i} \rangle \geq c \log s$.
\end{itemize}
As we will shortly see every step can still be carried out and we thus obtain Theorem~\ref{2Dlittlewood}. One way to describe this process is to say we will employ the McGehee--Pigno--Smith method in one dimension and the Cohen--Davenport--Pichorides in the other. 

The Cohen--Davenport--Pichorides method is applicable when one considers Freiman isomorphisms. We will thus employ it in all three dimensions to prove Theorem~\ref{Freiman Isomorphism Invariance}. The details can be found in the two upcoming sections. 

We begin with a technical result that is the main building block of the two proofs.
\begin{lemma} \label{Technical CDP}
Let $R$ and $d$ be positive integers, $K$ a positive real number and $F:\mathbb{T}^d \mapsto \mathbb{C}$ be an integrable function. Suppose there are positive integers $n_1,\dots,n_{R}$ and a collection of integrable functions $\Phi_{n_1},\dots ,\Phi_{n_{R}}$ such that
\begin{itemize}
\item[(A)] $\|\Phi_{n_i}\|_{\infty}\leq 1$ for $1\leq i \leq R$. 
\item[(B)] $\langle\Phi_{n_i}, F \rangle \geq K$ for $1\leq i \leq R$. 
\item[(C)] Let $l$ be a positive integer. $\langle F\,,\,  \Phi_{n_{i_0}} \Phi_{n_{i_1}} \cdots \Phi_{n_{i_{l}}} \overline{ \Phi_{n_{i_{l+1}}} \cdots \Phi_{n_{i_{2l}}}} \rangle = 0$ for $1\leq i_j \leq R$ unless $n_{i_0}+n_{i_1} + \dots + n_{i_l} - n_{i_{l+1}}- \dots - n_{i_{2l}}=n_\nu $ for some $1\leq \nu \leq R$.
\end{itemize}
Then there is a test function $g$ such that
\begin{itemize}
\item[(i)] $\|g\|_{\infty}\leq 1$. 
\item[(ii)] $g$ is a linear combination of functions of the form $\Phi_{n_{i_0}} \cdots \Phi_{n_{i_{k}}} \overline{ \Phi_{n_{i_{k}}} \cdots \Phi_{n_{i_{2k}}}}$ for some $k\leq 2 \log R$. 
\item[(iii)] $\langle g, F \rangle \geq c\, K \frac{\log^{1/2} R}{\log\log ^{1/2} R}\,.$
\end{itemize}
In particular (i) and (iii) imply
$$ \|F\|_1 \geq c \,K \frac{\log^{1/2}R}{\log\log^{1/2}R}\;.$$
\end{lemma}
The reader can think of the $\Phi_n$ as exponentials in order to gain some intuition. We will need two lemmata. The first is Lemma~1 of \cite{Pichorides1974}.
\begin{lemma}[Pichorides]\label{Pichorides}
Let $t\geq 100$. Suppose the quantities $P$ and $Q$ satisfy $t+2P\geq 0$ and $P^2+Q^2\leq t^4/4$. Then $$\left| 1-\frac{1}{t}-\frac{P+iQ}{t^2}+\frac{(P+iQ)^2}{t^4}\right|+\frac{1}{4t^{3/2}} (t+2P)^{1/2}\leq 1\;.$$
\end{lemma}
The second is also a result Pichorides (Lemma~2 in \cite{Pichorides1974}) whose proof is essentially due to Davenport (cf. Lemma~3 in \cite{Davenport1960}). 
\begin{lemma}[Davenport--Pichorides]\label{Davenport}
Let $E$ and $S$ be sets of positive integers. For $p\in S$ let $N(p)$ to be the number of elements of $E$ that are greater than $p$. 

Let $t$ be a positive integer and suppose that $$t^4 \sum_{p\in S} N(p)  \leq |E| \;.$$ Then there exist $t$ integers $\{ m_1, \dots , m_t\}$ in $E$ such that $$p+(m_\alpha-m_\beta)+(m_\gamma-m_\delta) \notin E$$ for all $p\in S$ and $1\leq \alpha\leq \beta \leq t,~1\leq \gamma <\delta\leq t$. 

Furthermore $m_\alpha=n_{q(\alpha)}$, where $q(\alpha)\leq \alpha^4 \sum_{p\in S} N(p)\,.$ 
\end{lemma}

We now turn to proving Lemma~\ref{Technical CDP}.
\begin{proof}[Proof of Lemma~\ref{Technical CDP}]
The proof is based on iteration. We will construct functions
$g_1,g_2,\ldots$ that satisfy (i) and modified versions of (ii) and
(iii):
\begin{itemize}
\item[($ii^\prime$)] $g_i$ is a linear combination of functions of
the form $\Phi_{n_{i_0}} \cdots \Phi_{n_{i_{k}}} \overline{
\Phi_{n_{i_{k+1}}} \cdots \Phi_{n_{i_{2k}}}}$ for some $k\leq 2
i$.
\item[($iii^\prime$)] $\langle g_i, F \rangle \geq  K (4
t^{1/2})^{-1} \sum_{n=0}^{i-1} (1-1/t)^n $ for some $t\geq 100$ to be
chosen later.
\end{itemize}
We set $g_1 = \Phi_{n_1}$, which satisfies (i),($ii^\prime$) and
($iii^\prime$) as the sum is empty. We now inductively define
\begin{eqnarray*}
g_{i+1}(x) &=& g_i(x)\left(1-\frac{1}{t}\right) \\ & & ~~~~~~~ -
g_i(x) \left(\frac{1}{t^2} \sum_{1\leq i < j\leq t}\Phi_{m_i}(x)
\overline{\Phi_{m_j}(x)} - \frac{1}{t^4} \left(\sum_{1\leq i <
j\leq t}\Phi_{m_i}(x) \overline{\Phi_{m_j}(x)}\right)^2 \right)
\\  & & ~~~~~~~~~~~~~~~~~~~~ + \frac{1}{4t^{3/2}} \sum_{1\leq i
\leq t}\Phi_{m_i}(x)
\end{eqnarray*}
for some $m_1,\dots,m_t$ carefully chosen from $\{n_1,\dots,n_R\}$
in such a way that the inner product of the middle part with $F$
is zero. For the time being we assume this can be done. We need to
check that $g_{i+1}$ satisfies (i),($ii^\prime$) and ($iii^\prime$).

For (i) we apply  Lemma~\ref{Pichorides}. For any $v$ set
$$P+iQ = \sum_{1\leq i<j \leq t} \Phi_{m_i}(v) \overline{\Phi_{m_j}(v)}\;.$$
We observe that
$$P^2+Q^2=|P+iQ|^2\leq \left(\sum_{1\leq i<j \leq t} |\Phi_{m_i}(v)|
|\Phi_{m_j}(v)|\right)^2 \leq \left(t(t-1)/2\right)^2<t^4/4$$ and
that $$0\leq \left|\sum_{i=1}^t \Phi_{m_i}(v)\right|^2=\sum_{i=1}^t |\Phi_{m_i}(v)|^2 +
2P \leq t+2P\;.$$ The conditions of Lemma~\ref{Pichorides} are
satisfied and so 
\begin{eqnarray*}
 |g_{i+1}(v)|  &\leq&  |g_i(v)| \left| 1-\frac{1}{t}-\frac{P+iQ}{t^2}+\frac{(P+iQ)^2}{t^4}\right|+\frac{1}{4t^{3/2}} (t+2P)^{1/2} \\
		    &\leq&  \left| 1-\frac{1}{t}-\frac{P+iQ}{t^2}+\frac{(P+iQ)^2}{t^4}\right|+\frac{1}{4t^{3/2}} (t+2P)^{1/2} \\ 
		    &\leq&  1\;.
\end{eqnarray*}
The last inequality coming from Lemma~\ref{Pichorides}. Thus $g_{i+1}$ satisfies (i). $g_{i+1}$ by definition satisfies ($ii^\prime$) and so we are left with
($iii^\prime$).

It follows from our assumption on the middle part of $g_{i+1}$
that $$\langle g_{i+1},F\rangle = \left(1-\frac{1}{t}\right)\langle
g_{i},F\rangle + \frac{1}{4 t^{3/2}} \sum_{i=1}^t \langle
\Phi_{m_i},F\rangle \geq \frac{K}{4 t^{1/2}} \sum_{n=0}^{i}
\left(1-\frac{1}{t}\right)^n.$$ Once $n$ becomes considerably bigger than $t$ the terms
$(1-1/t)^n \leq \exp(-n/t)$ become exponentially small and so add
very little to the sum. We therefore iterate the process only $t$ times and set $g=g_t$. It follows that the $k$ appearing in $(ii)$ can be taken to be $2 t$.
\begin{eqnarray}\label{Technical CDP - Lower Bound}
\langle g,F\rangle \geq \frac{K }{4 t^{1/2}} \sum_{n=1}^t
\left(1-\frac{1}{t}\right)^n \geq c K t^{1/2}
\end{eqnarray}
subject only to being able to repeat the iteration $t$ times. 

Our final task then becomes to prove that the $m_i$ can indeed be chosen $t$ times and get the largest possible value for $t$. This will be done by applying Lemma~\ref{Davenport}.

We start by labelling $m_1^{(i)},\dots,m_t^{(i)}$ the elements of
$\{n_1,\dots,n_R\}$ chosen in the $i$th iteration and recursively
define the following sets:
$$S_1=\{n_1\},~~S_{i+1}=S_i \cup T_i \cup U_i$$
where\\
$T_i=\{m_1^{(i)},\dots,m_t^{(i)}\}$ and\\
$U_i=\{p+(m_\alpha^{(i)}-m_\beta^{(i)})+(m_\gamma^{(i)}-m_\delta^{(i)}) :  p\in
S_i,\;1\leq \alpha \leq \beta \leq t,\; 1\leq \gamma < \delta \leq t\}.$

Let $E=\{n_1,\dots, n_R\}$. It follows from condition (C) that the middle part of $\langle
g_{i},F \rangle$ is zero provided that
$p+(m_\alpha^{(i)}-m_\beta^{(i)})+(m_\gamma^{(i)}-m_\delta^{(i)}) \notin E$ for all
$p \in S_{i-1}$, $1\leq \alpha \leq \beta \leq t$ and $1\leq \gamma < \delta \leq t$.

Applying Lemma~\ref{Davenport} with $S=S_{i-1}$ we see that the $m_j^{(i)}$ can be chosen provided that
$$t^4 \sum_{p\in S_{i-1}} N(p) \leq R\;.$$
The sum in the left hand side is estimated using the final conclusion of Lemma~\ref{Davenport}.
\begin{eqnarray*}
\sum_{p\in S_{i}} N(p) & = & \sum_{p\in S_{i-1}} N(p) + \sum_{p\in
                            T_{i-1}} N(p) + \sum_{p\in U_{i-1}} N(p)\\
                            & = & \sum_{p\in S_{i-1}} N(p) + \sum_{\alpha=1}^{t} N(m_\alpha^{(i)}) +
                            \sum_{p,\alpha,\beta,\gamma,\delta} N(p+(m_\alpha^{(i)}-m_\beta^{(i)})+(m_\gamma^{(i)}-m_\delta^{(i)}))\\
                            & \leq & \sum_{p\in S_{i-1}} N(p) + \sum_{\alpha=1}^{t} \alpha^4\left( \sum_{p\in S_{i-1}} N(p)\right) 
                            + t^4 \sum_{p \in S_{i-1}} N(p)\\
                            & \leq & t^5 \sum_{p \in S_{i-1}} N(p) \;.
\end{eqnarray*}
We used the fact that $m_\alpha^{(i)}-m_\beta^{(i)}+m_\gamma^{(i)}-m_\delta^{(i)}>0$ so that
$N(p+m_\alpha^{(i)}-m_\beta^{(i)}+m_\gamma^{(i)}-m_\delta^{(i)}) \leq N(p)$.

Observe that $\sum_{p\in S_1}N(p) = 1$. It follows by induction that $$ \sum_{p\in S_i} N(p)\leq t^{5i} \;.$$
The iteration is thus possible for $t$ steps when $t^{5t} \leq R$. So we take $$t= \left\lfloor \frac{\log R}{10\log\log R} \right\rfloor\;.$$  Substituting this value of $t$ in \eqref{Technical CDP - Lower Bound} gives conclusion $(iii)$. Conclusion $(ii)$  has been shown to hold for $k=2t\leq 2\log R$ and so has conclusion $(i)$. 
\end{proof}

\section[Towards a 2-dimensional Littlewood conjecture]{Towards a 2-dimensional Littlewood conjecture}

\label{2 dim Littlewood}

We now prove Theorem~\ref{2Dlittlewood}. Loosely speaking the
first dimension will be used to construct the $\Phi_n$ and the
second to combine them and produce a better test function. 
\begin{proof}[Proof of Theorem~\ref{2Dlittlewood}]
We apply Lemma~\ref{Technical CDP} to $F=F_A$. We take $e_1,e_2$
to be the standard basis of $\mathbb{Z}^2$ and translate $A$ if
necessary so that the coordinates of all its points are positive
integers. We let $A_1,\dots,A_R$ be the rows of $A$ and $n_i = A_i \cdot e_2$ for $1\leq i\leq R$. 

We set $\Phi_{n_i}$ to be the McGehee--Pigno--Smith test function for
$F_{A_i}$. By this we mean a function whose Fourier transform is
supported on $\{u\in \mathbb{Z}^2 : u\cdot e_2=n_i\}$ and which satisfies
$\|\Phi_{n_i}\|_\infty \leq 1$ and $\langle
F_{A_i},\Phi_{n_i}\rangle\geq c \log|A_i|\geq c \log s$. 

Hence the $\Phi_{n_i}$ satisfy conditions (A) and (B) for $K \geq c\log s$. Condition (C) is also satisfied as we see by examining the support of the Fourier transform of $\Phi_{n_{i_{0}}} \Phi_{n_{i_{1}}} \cdots \Phi_{n_{i_{l}}} \overline{ \Phi_{n_{i_{l+1}}} \cdots \Phi_{n_{i_{2l}}} }$: it lies on the line $\{ u\in \mathbb{Z}^2 : u\cdot e_2 = n_{i_0}+\dots+n_{i_l}-n_{i_{l+1}}-\dots-n_{i_{2l}} \}$. In particular $$\langle F_A \, ,\, \Phi_{n_{i_{0}}} \Phi_{n_{i_{1}}} \cdots \Phi_{n_{i_{l}}} \overline{ \Phi_{n_{i_{l+1}}} \cdots \Phi_{n_{i_{2l}}} } \rangle = 0$$ unless $n_{i_0}+\dots+n_{i_l}-n_{i_{l+1}}-\dots-n_{i_{2l}} = n_\nu$ for some $1\leq \nu \leq R$. The theorem follows from he final conclusion of Lemma~\ref{Technical CDP} by observing that $R\geq r$.
\end{proof}
We can of course take $\Phi_{n_i}$ to be the test function that
satisfies $\langle \Phi_{n_i},F_A\rangle = \|F_{A_i}\|_1$. Its Fourier transform is still supported on $\{u\in\mathbb{Z}^2 : u\cdot e_2 = n_i\}$ and hence everything we
did above can be repeated to yield the following.
\begin{theorem}
Let $A\subset \mathbb{Z}^2$ be finite. Suppose that $A$ consists of at least $r$ rows of size at least $s$. Then $$\|F_A\|_1 \geq c\, \mu(s) \left(\frac{\log r}{\log\log r}\right)^{1/2}$$ with $$ \mu(s)=\min\|F_{A_i}\|_1\;,$$where $A_1,A_2,\ldots$ are the rows of $A$.
\end{theorem}
%To deduce Corollary~\ref{Random} from Theorem~\ref{2Dlittlewood} we have to show that with high probability all the rows of $A$ have at least $N/4$ elements. This can be done by applying a famous inequality of S. Bernstein \cite{Bernstein1924}.
%\begin{lemma}[Bernstein] \label{Bernstein}
%Let $X_1,\dots,X_N$ be independent random variables with $\mathbb{P}(X_i=1)=\mathbb{P}(X_i=-1)=1/2$ for all $i$. Then for all $0<\eta<1$, $$\mathbb{P}\left(\left |  \sum_{i=1}^N X_i\right| \geq \eta N \right) \leq C \exp\left(- c N \eta^2 \right)\;.$$ 
%\end{lemma}
%\begin{proof}[Proof of Corollary~\ref{Random}]
%Let $A_1,\dots,A_N$ be the rows of $A$, where $A_j= \{v \in A : v\cdot e_2 = j\}$. Let $X_i^{j}$ be the random variable that takes value 1 or -1 according to whether $(i,j)\in A$. By Lemma~\ref{Bernstein} we have that for all $j$, $$ \mathbb{P}(\textrm{$A_j$ has at most $N/4$ elements}) \leq \mathbb{P}\left(\left |  \sum_{i=1}^N X_i^j\right| \geq  N/2 \right)\leq \exp(- c N)\;. $$ Thus  $$ \mathbb{P}(\textrm{$A_j$ has at least $N/4$ elements for all $1\leq j \leq N$}) \geq 1- N \exp(- c N)\;. $$ This probability tends to 1 as $N$ gets arbitrarily large. The lower bound on $\|F_A\|_1$ follows from Theorem~\ref{2Dlittlewood}.  
%\end{proof}

\section[Multidimensional sets in Z]{Multidimensional sets in $\mathbb{Z}$}

\label{Freiman Isomorphisms}

We repeat the same process to prove Theorem~\ref{Freiman
Isomorphism Invariance}. We can no longer use the
McGehee--Pigno--Smith test functions as their support is both very
large and very difficult to analyse. It is furthermore unlikely that condition (C) in Lemma~\ref{Technical CDP} holds. Instead we use the
Cohen--Davenport--Pichorides test functions, which are Freiman
isomorphism friendly because of conclusion (ii) in Lemma~\ref{Technical CDP}. In what follows for a set of integers $S$ and a positive integer $\alpha$ we write $$\alpha S = \{s_1+\dots+s_\alpha : s_i\in S\} \;.$$ 
\begin{proof}[Proof of Theorem~\ref{Freiman Isomorphism Invariance}]
Translate $A$ if necessary so that all three coordinates of its elements are positive. Let $\theta$ be the Freiman isomorphism between $A$ and $B$ and
$e_1,e_2,e_3$ the standard basis of $\mathbb{Z}^3$. Suppose that
$A_{1}, A_2,\dots$ are the planar slices of $A$. For any $i$ let $a_i$ be the integer such that $A_i \subset \{u\in\mathbb{Z}^3: u\cdot e_3 = a_i\}$. Each $A_i$ consists of at least $r$ rows $A_{i1},A_{i2}\dots$ of size at least $s$. Let $a_{ij}$ be the integer such that  $A_{ij}\subset \{u\in\mathbb{Z}^3: u\cdot e_3=a_i,\, u\cdot e_2=b^i_j\}$. 

We construct a test function for $F_B=F_{\theta(A)}$ by three successive applications of Lemma~\ref{Technical CDP}.

We begin by applying Lemma~\ref{Technical CDP} to get a test
function for $F_{\theta(A_{ij})}$ for all pairs of indices $\{i,j\}$ for which $A_{ij}$ is non-empty. Let $b_{ij}^{(1)},
b_{ij}^{(2)}, \ldots$ be the elements of $\theta(A_{ij})$. We set $n_l=b_{ij}^{(l)}$ and $\Phi_{n_l} =
e(b_{ij}^{(l)})$ in Lemma~\ref{Technical CDP}. The $\Phi_{n_l}$ satisfy conditions (A), (B)
with $K=1$ and (C). Applying Lemma~\ref{Technical CDP} we get a
test function $ f_{ij}$ which satisfies $$\langle
F_{\theta(A_{ij})}, f_{ij}\rangle\geq c \left(\frac{\log s}{\log\log s}\right)^{1/2}$$
and $\|f_{ij}\|_\infty\leq 1$. Next we observe that the support of
$\widehat{f_{ij}}$ lies in $(\alpha+1)\theta(A_{ij})-\alpha\,\theta(A_{ij})$
for some $\alpha  \leq 2 \log s$. In particular it does not intersect
$\theta(A \setminus A_{ij})$, for if $\theta(u)=
\theta(u_0)+\dots+\theta(u_\alpha)-\theta(u_{\alpha+1})-\dots -
\theta(u_{2\alpha})$ for some $u\in A\setminus A_{ij}$ and $u_0,\dots,u_{2\alpha}\in
A_{ij}$, then $u = u_0+\dots+u_\alpha-u_{\alpha+1}-\dots -u_{2\alpha}$ as $\theta$ is a Freiman isomorphism of degree $k$ and $\alpha\leq
k$. This is impossible as the right hand side is supported on the
line $\{u\in\mathbb{Z}^3:  u\cdot e_3=a_i,\, u\cdot e_2=b^i_j\}$, while the left hand
is not. Hence $$\langle F_{\theta(A_i)} , f_{ij}, \rangle = \langle
F_{\theta(A_{ij})},f_{ij}\rangle\geq c \left(\frac{\log s}{\log\log s}\right)^{1/2}\;.$$

Next we combine the $f_{ij}$ to get a test function for
$F_{\theta(A_i)}$. We set $n_j=b_j^{i}$ and $\Phi_{n_j}=f_{ij}$ in Lemma~\ref{Technical CDP}.
The $f_{ij}$ satisfy condition (A) and, as we saw above, (B) with $K\geq c \log ^{1/2-\varepsilon}s$. To
check condition (C) note that the Fourier transform of
$f_{ij_{0}}f_{ij_{1}}f_{ij_{2}}\cdots
f_{ij_{l}}\overline{f_{ij_{l+1}}f_{ij_{l+2}} \cdots f_{ij_{2l}}}$
is supported on 
\begin{eqnarray*}
(\alpha+1) \theta(A_{ij_{0}})- \alpha \theta(A_{ij_{0}}) + (\alpha+1) \theta (A_{ij_{1}}) -\alpha\theta (A_{ij_{1}}) + \dots + (\alpha+1) \theta(A_{ij_{l}}) - \\
\alpha \theta(A_{ij_{l}})- (\alpha+1) \theta({A_{ij_{l+1}})+\alpha\theta({A_{ij_{l+1}})-\dots - (\alpha+1) \theta(A_{ij_{2l}}})+\alpha \theta(A_{ij_{2l}}})
\end{eqnarray*}
for $\alpha \leq 2\log s$.
Thus the inner product with $F_{\theta(A_i)}$ is zero unless
$\theta(A_i)$ intersects the above sum-difference set. Note that
$l\leq 2 \log r$ and that $\theta$ is a Freiman isomorphism of
sufficiently large degree for this to happen only when the sum $
b_{j_0}^i+b_{j_1}^i+\dots+b_{j_l}^i - b_{j_{l+1}}^i-\dots -
b_{j_{2l}}^i$ equals  $b_j^i$ for some $j$. 

By Lemma~\ref{Technical CDP} we get a test function $f_i$ that satisfies $\|f_i\|_\infty\leq 1$ and $$\langle
F_{\theta(A_i)},f_{i}\rangle\geq c \left(\frac{\log s\,\log r}{\log\log s\, \log\log r}\right)^{1/2}\;.$$ The support
of $\widehat{f_i}$ lies in $(\gamma+1)\theta(A_i)-\gamma\theta(A_i)$ for
some $\gamma\leq 12 \log r \log s$: the support of $\widehat{f_{ij}}$
lies in $(\alpha+1)\theta(A_i)-\alpha\theta(A_i)$ for $\alpha\leq 2 \log s$ and
we have to consider expressions of the form
$f_{ij_{0}}f_{ij_{1}}f_{ij_{2}}\cdots
f_{ij_{\beta}}\overline{f_{ij_{\beta+1}}f_{ij_{\beta+2}} \cdots f_{ij_{2\beta}}}$
for $\beta\leq 2 \log r$ and so $\gamma$ can be taken to be
$(\alpha+1)\,\beta+\alpha\,(\beta+1)=2\alpha\beta+\alpha+\beta \leq 12 \log r \log s$. Thus the support
of $\widehat{f_i}$ does not intersect $\theta(A\setminus A_i)$, for if $\theta(u)=
\theta(u_0)+\dots+\theta(u_\gamma)-\theta(u_{\gamma+1})-\dots -
\theta(u_{2\gamma})$ for some $u\in A\setminus A_i$, $u_l\in A_{i}$
and $\gamma\leq 12 \log s\log r$, then, as $\theta$ is a Freiman
isomorphism of degree $k\geq \gamma$, $u$ would have to equal
$u_0+\dots+u_\gamma-u_{\gamma+1}-\dots -u_{2\gamma}$. This is impossible as the
right hand side lies on the plane $\{u \in \mathbb{Z}^3 : u\cdot e_3=a_i\}$, while the
left hand does not. Hence $$\langle F_{\theta(A)} , f_{i}, \rangle
= \langle F_{\theta(A_{i})},f_{i}\rangle\geq c \left(\frac{\log s\,\log r}{\log\log s\, \log\log r}\right)^{1/2}\;.$$

Finally we combine the $f_i$ to get a test function for
$F_{\theta(A)}$. We let $n_i=a_{i}$ and $\Phi_{n_i}=f_{i}$. The $f_{i}$ satisfy conditions (A) and, 
as we saw above, (B) with $K\geq c(\log s \,\log
r)^{1/2-\varepsilon}$ in the statement of Lemma~\ref{Technical CDP}. To check condition (C) note that the Fourier
transform of $f_{i_{0}}f_{i_{1}}f_{i_{2}}\cdots
f_{i_{l}}\overline{f_{i_{l+1}}f_{i_{l+2}} \cdots f_{i_{2l}}}$ is
supported on $(\gamma+1) \theta(A_{i_{0}})-\gamma \theta(A_{i_{0}}) + (\gamma+1)
\theta (A_{i_{1}}) - \gamma \theta (A_{i_{1}}) +\dots +(\gamma+1)
\theta(A_{i_{l}}) - \gamma \theta(A_{i_{l}})- (\gamma+1)
\theta(A_{i_{l+1}})+\gamma\theta(A_{i_{l+1}})-\dots -(\gamma+1)
\theta(A_{i_{2l}})+\gamma\theta(A_{i_{2l}})$ for $\gamma\leq 12 \log s\log
r$. The inner product with $F_{\theta(A)}$ is zero unless
$\theta(A)$ intersects the above sum-difference set, which is a
subset of $(\delta+1)\theta(A) - \delta\theta(A)$ for $\delta = 2l\gamma+l+\gamma\leq 62
\log p \log r \log s$. $\theta$ is a Freiman isomorphism of degree
$k\geq \delta $, so this happens only if $ a_{i_0}+a_{i_1}+\dots+a_{i_l} - a_{i_{l+1}}-\dots -a_{i_{2l}}$
equals  $a_i$ for some $i$. By Lemma \ref{Technical CDP} we get $$\|F_{\theta(A)}\|_1 \geq c \left(\frac{\log s\,\log r\,\log p}{\log\log s\, \log\log r \,\log\log p}\right)^{1/2} \;.\qedhere$$
\end{proof}
\begin{remark}
One can extend this result to higher dimensions. 
\end{remark}

\section[Additive structure when the L^1-norm is small]{Additive structure when $\|F_A\|_1$ is small}

\label{Additive Structure when L1 Norm Small}

In this final section we discuss the following question. Suppose $\|F_A\|\leq C \log|A|$ for $A\subset \mathbb{Z}$. Is there a particular structure $A$ must have? We suggest a plausible structure and compare it with that implied by Theorem~\ref{Freiman Isomorphism Invariance}.

Determining the precise value of $\|F_A\|_1$ for a given $A$ is hard. The Cauchy-Schwarz inequality shows that the $L^1$-norm is certainly bounded above by the $L^2$-norm,
$\|F_A\|_2 = |A|^{1/2}$. This order of magnitude is attained when $A$ is the lacunary sequence $\{2^i : 1\leq i\leq N\}$. By an averaging argument one gets much
denser random subsets of $\{1,2,\dots,N\}$ with $\|F_A\|_1 \geq c
N^{1/2}$. In general sets with random like properties are expected
to give rise to exponential sums with large $L^1$-norm. For example,  if $A$ is the set of the first
$N$ primes, then $\|F_A\|_1 \geq N^{1/2 - \varepsilon}$ for all $\varepsilon>0$
\cite{Vaughan1988} and if $A$ is the intersection of the support
of the M\"{o}bius function with $\{1,2,\dots,N\}$, then $\|F_A\|_1
\geq N^{1/8 - \varepsilon}$ \cite{Balog-Ruzsa1999}.

At the other end of the spectrum we have structured sets. If $A$ is the union of $k$ arithmetic progressions, then by the triangle inequality $\|F_A\|_1 \leq C\,k \log|A|$. Furthermore, if $A$ is a $d$-dimensional arithmetic progression $$\{c+x_1q_1+\dots+ x_d q_d: 0\leq x_i\leq N\,\mbox{for $1\leq i \leq d$}\} ~\mbox{for $c, q_i\in \mathbb{Z}$ for $1\leq i \leq d$}\;,$$ then $\|F_A\|_1 \leq (C \log|A|)^d$. 

Note however that not the whole of $A$ needs to be structured. We can for example remove a subset $X$ with $C \log^2N$ elements from $\{1,2,\dots,N\}$ and still have $$\|F_A\|_1 \leq  \|F_{\{1,\dots,N\}}\|_1 + \|F_B\|_1 \leq C \log N + \|F_B\|_2 = C \log N + C\log N \leq C \log |A|\;.$$ One can instead add a much larger set $X$. For example $X$ can be a 2-dimensional arithmetic progression disjoint from $\{1,\dots,N\}$. If $X$ is Freiman 2-isomorphic to $\{1,\dots,L\}\times\{1,\dots,L\}$, where $L=\exp(\log^{1/2}N)$, then $\|F_X\|_1\leq C \log N$. Thus $\|F_{\{1,\dots,N\}\cup X}\|_1 \leq C \log|A|$. 

Establishing a concrete relation between $\|F_A\|_1$ and the additive structure of $A$ has not been possible so far. Even the simplest inverse theorem for sets $A$ where $\|F_A\|_1$ is close to being minimal has been elusive. The following question arose in conversations with B.J. Green and is in accordance with a theorem of Green and T. Sanders on idempotent measures \cite{Green-Sanders2008}. 
\begin{question}\label{ConjectureGreen}. %To state it les us write $1_S$ for the characteristic function of a set of integers $S$.
Does there exists an absolute constant $1/2 \leq \eta <1$ and a function $g: \mathbb{R}^+\mapsto\mathbb{R}^+$ with the following property. Let $A\subset \mathbb{Z}$ be a finite set and $K$ a positive constant. Suppose $\|F_A\|_1\leq K \log|A|$. Then there exists a set $X\subset \mathbb{Z}$ of size at most $\exp(\|F_A\|_1^ \eta)$, $g(K)$ arithmetic progressions $P_1,\dots,P_{g(K)}$ and $\varepsilon_1,\dots,\varepsilon_{g(K)}\in \{+1,-1\}$ such that $$F_A= F_X + \sum_{i=1}^{g(K)} \varepsilon_i F_{P_i} \;.$$
\end{question} 
The range of $ \eta$ comes from the example discussed above and Theorem~\ref{MPSG}. Taking $A$ to be a 2-dimensional arithmetic progression Freiman 2-isomorphic to $\{1,\dots,N\}\times\{1,\dots,N\}$ suggests that $g(K)$ has to be exponential in $K$.

The results in this paper point to a slightly different direction. We have established that no sufficiently large set of integers $A$ whose exponential sum has $L^1$-norm at most $C\log|A|$ can be Freiman isomorphic to a genuinely three dimensional set in $\mathbb{Z}^3$. This puts a constraint on sets where $\|F_A\|_1$ is close to being minimal. Unfortunately it is not the case that such sets mainly consist of few long arithmetic progressions and a small set. The notion of dimensionality we have relied on is too restrictive to lead to such a conclusion. 

Take for example the lacunary sequence $A=\{x_i=2^i : 1\leq i\leq N\}$. Its elements satisfy the recurrence relation $x_{i+1}=x_i+2(x_i-x_{i-1})$. It follows that its image under a Freiman isomorphism $\theta$ of degree 3 also satisfies this relation. The $y$-coordinate of the elements of $\theta(A)$ is either constant (when $\theta(x_1)\cdot e_2 =\theta(x_2)\cdot e_2$) or distinct for all $i$. In other words either $\theta(A)$ is contained in a single row or it consists of $|A|$ singleton rows. In either case $\theta(A)$ is not a genuinely 3-dimensional set. Yet any subset $Y\subset A$ cannot be decomposed in fewer than $|Y|/2$ arithmetic progressions as $A$ contains at most two consecutive elements of any arithmetic progression.  

Lacunary sequences are very sparse, but the situation doesn't change when we consider dense sets as the following example demonstrates.

Let $L$ be a large integer and $P$ the first prime such that
$$\sum_{p \in \mathcal{P} }p^{-1} \geq 1/2$$
where $\mathcal{P}$ is the set of primes between $L$ and $P$. Now let
$$ N = \prod_{p\in \mathcal{P}} p$$
and
$$A = \bigcup_{p\in \mathcal{P}} A_p \;,$$
where $A_p$ consists of all numbers in $\{1,\dots, N\}$ that are congruent to $1 \bmod p$.

$A$ has large density in $\{1,\dots,N\}$. To check this observe that
$$ |A| = \left| \bigcup_{p\in \mathcal{P}} \left(A_p / \bigcup_{q\neq p} A_q\right)\right| = \sum_{p\in \mathcal{P}} \left|A_p / \bigcup_{q\neq p} A_q\right|\;.$$
We know that $|A_p|=N/p$ and $|A_p \cap A_q| = N/pq$. Hence
$$\left|A_p / \bigcup_{q\neq p} A_q\right|\geq \frac{N}{p} \,
\left(1-\sum_{q\in \mathcal{P}\setminus \{p\}} q^{-1}\right)\geq
\frac{N}{2p}\;.$$ Which in turn implies that
$$|A|\geq \frac{N}{2} \sum_{p\in \mathcal{P}} p^{-1}\geq  N/4\;.$$
Next we consider the image of $A$ under a Freiman isomorphism of degree two.
Freiman isomorphisms map arithmetic progressions in $\mathbb{Z}$
into lines in $\mathbb{Z}^3$ and hence $\theta(A)$ must be
supported on a collection of lines $\{\theta(A_p) : p\in
\mathcal{P}\}$. For every pair of indices $p\neq q$, $\theta(A_p) \cap \theta(A_q) = N/pq >
2$ and so the two lines must in fact be identical. Thus the image
of $A$ under any Freiman isomorphism lies in a single line in
$\mathbb{Z}^3$. As a consequence $\theta(A)$ either lies in a single row or in $|A|$
different rows. 

\bibliography{all}

\begin{thebibliography}{10}

\bibitem{Balog-Ruzsa1999}
A.~Balog and {I.\,Z.} Ruzsa.
\newblock A new lower bound for the ${L^1}$ mean of exponential sums with the
  {M}{\"{o}}bius function.
\newblock {\em Bull. Lond. Math. Soc.}, 31:415--418, 1999.

\bibitem{Cohen1960}
{P.\,J.} Cohen.
\newblock On a conjecture of {L}ittlewood and idempotent measures.
\newblock {\em Amer. J. Math}, 82:191--212, 1960.

\bibitem{Davenport1960}
H.~Davenport.
\newblock On a theorem of {P}.\,{J}. {C}ohen.
\newblock {\em Mathematika}, 7:93--97, 1960.

\bibitem{Green-Sanders2008}
{B.\,J.} Green and T.~Sanders.
\newblock A quantitative version of the idempotent theorem in harmonic
  analysis.
\newblock {\em Ann. of Math.}, 168(3):1025--1054, 2008.

\bibitem{Littlewood1948}
{G.\,H.} Hardy and {J.\,E.} Littlewood.
\newblock A new proof of a theorem on rearrangements.
\newblock {\em J. London Math. Soc.}, 23:163--168, 1948.

\bibitem{Konyagin1981}
{S.\,V.} Konyagin.
\newblock On the {L}ittlewood problem.
\newblock {\em Izv. Akad. Nauk SSSR Ser. Mat.}, 45:243--265, 1981.

\bibitem{Lorch1954}
L.~Lorch.
\newblock The principal term in the asymptotic expansion of the {L}ebesgue
  constants.
\newblock {\em Amer. Math. Monthly}, 61:245--249, 1954.

\bibitem{MPS1981}
{O.\,C.} McGehee, L.~Pigno, and B.~Smith.
\newblock Hardy's inequality and the ${L^1}$ norm of exponential sums.
\newblock {\em Annal. of Math.}, 113(3):613--618, 1981.

\bibitem{Pichorides1974}
{S.\,K.} Pichorides.
\newblock A lower bound for the ${L^1}$ norm of exponential sums.
\newblock {\em Mathematika}, 21:155--159, 1974.

\bibitem{Vaughan1988}
{R.\,C.} Vaughan.
\newblock The ${L^1}$ mean of exponential sums over primes.
\newblock {\em Bull. Lond. Math. Soc.}, 20(2):121--123, 1988.

\end{thebibliography}

$\hspace{12pt}$ \textit{Email address}: giorgis@cantab.net

\end{document}